\documentclass[english]{cccconf}
\usepackage[comma,numbers,square,sort&compress]{natbib}
\usepackage{amsmath}
\usepackage{algorithmicx}
\usepackage{algpseudocode}
\usepackage{booktabs}
\usepackage{subfig}
\usepackage{fancyhdr}
\begin{document}
\title{State transition algorithm for traveling salesman problem}
\author{YANG Chunhua\aref{csu},
        TANG Xiaolin\aref{csu},
        ZHOU Xiaojun\aref{csu,ub},
        GUI Weihua\aref{csu}}



\affiliation[csu]{School of Information Science and Engineering,
        Central South University, Changsha 410083, P.~R.~China}
\affiliation[ub]{School of Science, Information Technology and Engineering,
         University of Ballarat, Victoria 3353, Australia
        \email{tiezhongyu2010@gmail.com}}
\maketitle

\begin{abstract}
Discrete version of state transition algorithm is proposed in order to solve the traveling salesman problem. Three special operators for discrete optimization problem named swap, shift and symmetry transformations are presented. Convergence analysis and time complexity of the algorithm are also considered.
To make the algorithm simple and efficient, no parameter adjusting is suggested in current version. Experiments are carried out to test the performance of the strategy, and comparisons with simulated annealing and ant colony optimization have demonstrated the effectiveness of the proposed algorithm. The results also show that the discrete state transition algorithm consumes much less time and has better search ability than its counterparts, which indicates that state transition algorithm is with strong adaptability.
\end{abstract}


\footnotetext{YANG Chunhua's work is supported by the National Science Found
for Distinguished Young Scholars of China (Grant No.
61025015), and ZHOU Xiaojun is supported by China Scholarship Council, who is the corresponding author of this paper.}

\section{Introduction}
As one of the most significant combinatorial optimization problems, traveling salesman problem (TSP) has been paid greatly attention and been extensively studied, due to its importance in manufacturing, distribution management and scheduling\cite{1,2}. Traveling salesman problem can
be described as: Given a set of \textit{n} nodes and distances for each pair of nodes, find a roundtrip of minimal total length visiting each node exactly once, and according to whether the distance from node \textit{i} to node \textit{j} is the same as from node \textit{j} to node \textit{i}, the TSP can be symmetric or asymmetric\cite{3}. For a \textit{n} nodes asymmetric TSP, there exists (\textit{n}-1)! possible solutions; while for the corresponding symmetric TSP, (\textit{n}-1)!/2 solutions are possible. Because of its NP-hard property, emphasis has shifted from the aim of finding a global optimal solution to the goal of obtaining ``good solutions" in reasonable time and establishing the ``the degree of goodness"\cite{4,5}.\\
\indent Taking the impractical exhaustive search for TSP into consideration, heuristic algorithms are introduced to speed up the process of finding a satisfactory solution, of which, simulated annealing (SA), tabu search (TS), genetic algorithm (GA), ant colony optimization (ACO), particle swarm optimization (PSO) have found  widest applications in the the field\cite{6,7,8,9,10}. In terms of the concepts of state and state transition, a new heuristic search algorithm named state transition algorithm (STA) is proposed recently, and it exhibits good search performance in continuous space\cite{11,12}. In the initial version of state transition algorithm, a special transformation named general elementary matrix was also proposed to solve the discrete optimization problem. In this paper, we will continue to develop its detailed forms to promote deep study of STA for solving traveling salesman problem.\\
\indent The paper is organized as follows: In section 2, we review the unified form of the state transition algorithm, and establish the discrete version of STA. Then detailed forms and implementations of three special transformation operators are put forward in section 3, and we also give some strategies of parameters setting. In the next section, some experiments are performed to test the performance of the proposed algorithm. Conclusion is derived in the end.
\section{State transition algorithm in discrete version}
As stated in\cite{11}, a solution to a specific optimization problem can be described as a state, the thought of optimization algorithms can be treated as state transition, and the process to update solutions will become a state transition process. Without loss of generality, the unified form of state transition algorithm can be described as
\begin{eqnarray}
\left \{ \begin{array}{ll}
x_{k+1}= A_{k}x_{k} + B_{k}u_{k}\\
y_{k+1}= f(x_{k+1})
\end{array} \right.
\end{eqnarray}
where $x_{k} \in \Re^{n}$ stands for a state, corresponding to a solution of a optimization problem; $A_{k}$,
$B_{k} \in \Re^{n \times n}$ are state transition matrixes, which are usually transformation operators;
$u_{k} \in \Re^{n}$ is the function with variables $x_{k}$ and history states; $f$ is the cost function or evaluation function.\\
\indent As for traveling salesman problem, the cost $f$ is usually expressed as the function of a sequence, which corresponds to an ordered traveling route. That's to say, after a transformation using  $A_{k}$ or  $B_{k}$, a new state $x_{k+1}$ should be a sequence too. Due to the particularity of TSP, only a state transition matrix is considered, avoiding the complexity of ``adding" one sequence to another.\\
\indent According to the theory of Linear Algebra, we know that the special matrix must have only one position with value 1 in each column and each row. The special random matrix is called general elementary matrix, because it is derived from identity matrix and has the function of transforming a sequence into another one. For simplicity and specificity, $G_{k}$ is denoted in state transition algorithm for discrete optimization problem as follows
\begin{eqnarray}
\left \{ \begin{array}{ll}
x_{k+1}= G_{k}x_{k}\\
y_{k+1}= f(x_{k+1})
\end{array} \right.
\end{eqnarray}
where ${x}_k = [x_{1k},x_{ik},\cdots,x_{nk}]^T,x_{ik} \in \{1,2,\cdots,n \}$ be a sequence, the $G_{k}$ is general elementary matrix with only position value 1 in each column and each row.
For instance, $G_k$ has the following styles\cite{11}
\begin{tiny}
\begin{center}
$\left(
  \begin{array}{ccccc}
    1   & 0 & 0 & 0 & 0 \\
    0   & 1 & 0 & 0 & 0 \\
    0   & 0 & 0 & 0 & 1 \\
    0   & 0 & 0 & 1 & 0 \\
    0   & 0 & 1 & 0 & 0 \\
  \end{array}
\right)$
$\left(
  \begin{array}{ccccc}
    0   & 0 & 0 & 1 & 0 \\
    0   & 1 & 0 & 0 & 0 \\
    0   & 0 & 1 & 0 & 0 \\
    0   & 0 & 0 & 0 & 1 \\
    1   & 0 & 0 & 0 & 0 \\
  \end{array}
\right)$
$\left(
  \begin{array}{ccccc}
    0   & 0 & 0 & 1 & 0 \\
    0   & 1 & 0 & 0 & 0 \\
    1   & 0 & 0 & 0 & 0 \\
    0   & 0 & 0 & 0 & 1 \\
    0   & 0 & 1 & 0 & 0 \\
  \end{array}
\right)$
\end{center}
\end{tiny}
If an initial sequence is [1,2,3,4,5], then after the above general elementary matrices, it will becomes [1,2,5,4,3], [4,2,3,5,1], [4,2,1,5,3], respectively.
\section{State transition operators and corresponding adjusting strategies}
As can be seen from above, the general elementary matrix can have various types, and arbitrary types may have neither validity nor efficiency to discrete optimization problem. To make the state transition process controllable and efficient, three special transformation operators for traveling salesman problem are defined.\\
\noindent (1) Swap transformation
\begin{eqnarray}
x_{k+1}= G^{swap}_{k}(m_a) x_{k}
\end{eqnarray}
where $G^{swap}_{k} \in \Re^{n \times n}$ is called swap transformation matrix, $m_a$ is a constant called swap factor to control the maximum number of random positions to be exchanged, while the positions are random. For example, if $n=5, m_a=3$, the swap transformation matrix has the following styles
\begin{tiny}
\begin{center}
$\left(
  \begin{array}{ccccc}
     0  &   0  &   1  &   0   &  0\\
     0  &   1  &   0  &   0   &  0\\
     1  &   0  &   0  &   0   &  0\\
     0  &   0  &   0  &   1   &  0\\
     0  &   0  &  0   &   0   &  1\\
  \end{array}
\right)$
$\left(
  \begin{array}{ccccc}
     1  &   0  &   0  &   0   &  0\\
     0  &   0  &   1  &   0   &  0\\
     0  &   1  &   0  &   0   &  0\\
     0  &   0  &   0  &   1   &  0\\
     0  &   0  &   0  &   0   &  1\\
  \end{array}
\right)$
$\left(
  \begin{array}{ccccc}
     1  &   0   &  0   &  0  &   0\\
     0  &   0   &  0   &  0  &   1\\
     0  &   0   &  1   &  0  &   0\\
     0  &   1   &  0   &  0  &   0\\
     0  &   0   &  0   &  1  &   0\\
  \end{array}
\right)$
\end{center}
\end{tiny}
If an initial sequence is [1,2,3,4,5], then after the above swap transformation matrices, it will becomes [3,2,1,4,5], [1,3,2,4,5], [1,5,3,2,4], respectively.
It can be found that the 1\textit{th} and 3\textit{th} rows are exchanged, the 2\textit{th} and 3\textit{th} rows are exchanged, and the 2\textit{th}, 4\textit{th}, and 5\textit{th} rows are interchanged, respectively, which indicates that the maximum number of 3 positions are exchanged.\\
\noindent (2) Shift transformation\\
\begin{eqnarray}
x_{k+1}= G^{shift}_{k}(m_b) x_{k}
\end{eqnarray}
where $G^{shift}_{k} \in \Re^{n \times n}$ is called shift transformation matrix, $m_b$ is a constant called shift factor to control the maximum length of consecutive positions to be shifted. By the way, the selected position to be shifted after and positions to be shifted are chosen randomly. To make it more clearly, let $n=5, m_b=2$, then the shift transformation matrix has the following styles
\begin{tiny}
\begin{center}
$\left(
  \begin{array}{ccccc}
     1  &   0  &   0  &   0   &  0\\
     0  &   0  &   1  &   0   &  0\\
     0  &   1  &   0  &   0   &  0\\
     0  &   0  &   0  &   1   &  0\\
     0  &   0  &  0   &   0   &  1\\
  \end{array}
\right)$
$\left(
  \begin{array}{ccccc}
     1  &   0  &   0  &   0   &  0\\
     0  &   0  &   1  &   0   &  0\\
     0  &   0  &   0  &   1   &  0\\
     0  &   1  &   0  &   0   &  0\\
     0  &   0  &   0  &   0   &  1\\
  \end{array}
\right)$
$\left(
  \begin{array}{ccccc}
     1  &   0   &  0   &  0  &   0\\
     0  &   0   &  0   &  1  &   0\\
     0  &   1   &  0   &  0  &   0\\
     0  &   0   &  1   &  0  &   0\\
     0  &   0   &  0   &  0  &   1\\
  \end{array}
\right)$
\end{center}
\end{tiny}
If an initial sequence is [1,2,3,4,5], then after the above shift transformation matrices, it will becomes [1,3,2,4,5], [1,3,4,2,5], [1,4,2,3,5], respectively.
In the first two cases, the position to be shifted is \{2\}, and the positions to be shifted after are \{3\},\{4\}, while in the third case, the positions to be shifted is \{2,3\}, and the position to be shifted after are \{4\}. In the last case,
\{2,3\} are the maximum length of consecutive positions to be shifted.\\
\noindent (3) Symmetry transformation
\begin{eqnarray}
x_{k+1}= G^{sym}_{k}(m_c) x_{k}
\end{eqnarray}
where $G^{sym}_{k} \in \Re^{n \times n}$ is called symmetry transformation matrix, $m_c$ is a constant called symmetry factor to control the maximum length of subsequent positions as center.
By the way, the component before the subsequent positions and consecutive positions to be symmetrized are both created randomly. For instance, $n=5, m_c=1$, then the symmetry transformation matrix has the following styles
\begin{tiny}
\begin{center}
$\left(
  \begin{array}{ccccc}
     1  &   0  &   0  &   0   &  0\\
     0  &   1  &   0  &   0   &  0\\
     0  &   0  &   0  &   1   &  0\\
     0  &   0  &   1  &   0   &  0\\
     0  &   0  &   0  &   0   &  1\\
  \end{array}
\right)$
$\left(
  \begin{array}{ccccc}
     1  &   0  &   0  &   0   &  0\\
     0  &   0  &   0  &   0   &  1\\
     0  &   0  &   0  &   1   &  0\\
     0  &   0  &   1  &   0   &  0\\
     0  &   1  &   0  &   0   &  0\\
  \end{array}
\right)$
$\left(
  \begin{array}{ccccc}
     1  &   0  &   0  &   0   &  0\\
     0  &   1  &   0  &   0   &  0\\
     0  &   0  &   0  &   0   &  1\\
     0  &   0  &   0  &   1   &  0\\
     0  &   0  &   1  &   0   &  0\\
  \end{array}
\right)$
\end{center}
\end{tiny}
If an initial sequence is [1,2,3,4,5], then after the above symmetry transformation matrices, it will becomes [1,2,4,3,5], [1,5,4,3,2], [1,2,5,4,3], respectively.
All of the three cases, the component before the subsequent positions is \{3\}. In the first two cases, the subsequent position is \{$\emptyset$\}, while the
consecutive positions to be symmetrized are \{4\} and \{4,5\}, respectively. In the third case, the subsequent position(or the center) is \{4\}, while the consecutive position to be symmetrized is \{5\}. It is not difficult to find that the length of subsequent positions is 1 in the last case, which is the maximum length of subsequent positions as indicated by $m_c$.

\indent To accept a new solution, ``greedy criterion" is commonly adopted, in other words, only solution better than previous one is accepted. However, in SA, a worse solution is accepted probabilistically. Although the strategy in SA can jump out of stagnant state and it creates a big perturbation to previous solution, to guarantee the convergence of SA becomes difficult. In current version of state transition algorithm for discrete optimization problem, ``greedy criterion" is inherited. The core procedure of the STA for TSP can be outlined in pseudocode as follows
\begin{algorithmic}[1]
\Repeat
   \State [\textit{Best},\textit{fBest}] $\gets$ op\_swap(\textit{cities},\textit{Best},\textit{fBest},\textit{SE},\textit{n},$m_a$)
   \State [\textit{Best},\textit{fBest}] $\gets$ op\_shift(\textit{cities},\textit{Best},\textit{fBest},\textit{SE},\textit{n},$m_b$)
   \State [\textit{Best},\textit{fBest}] $\gets$ op\_symmetry(\textit{cities},\textit{Best},\textit{fBest},\textit{SE},\textit{n},$m_c$)
\Until{the specified termination criterion is met}
\end{algorithmic}
where, \textit{cities} is the information about the TSP; the \textit{SE} is the times of transformation, called search enforcement. As for more detailed explanations, op\_swap function of above pseudocode is also given in MATLAB scripts
\begin{verbatim}
for i=1:SE
    Tranf = swap_matrix(n,ma);
    State(i,:)= (Tranf*Best')';
end
[newBest,fGBest] = fitness(State,cities,SE);
if fGBest < fBest
    Best = newBest;
    fBest = fGBest;
end
function y = swap_matrix(n,ma)
y  = eye(n);
R  = randperm(n);
T  = R(1:ma);
S  = T(randperm(ma));
y(T,:) = y(S,:);
\end{verbatim}
\section{Theoretical analysis of discrete STA}
Next, we analyze the convergence performance, the global search ability, and time complexity of the discrete state transition algorithm.
\subsection{Convergence properties of discrete STA}
Firstly, we define the concept of convergence for discrete STA
\begin{eqnarray}
|f(x_k) - f(x^{*})| \leq \epsilon, \forall ~ k > N
\end{eqnarray}
while $x^{*}$ is the global minimum solution of a traveling salesman problem, $\epsilon$ is a small constant, and $N$ is a natural number. If $\epsilon > 0$, we can say that the algorithm converges to a $\epsilon-optimal~solution$; if $\epsilon = 0$, we can say that the algorithm converges to a \textit{global minimum}. However, if $x^{*}$ is a local minimum solution, then we can say that the algorithm converges to a $\epsilon-suboptimal~solution$ and a \textit{local minimum} for $\epsilon > 0$ and $\epsilon = 0$, respectively.\\
\textbf{Theorem 1} the discrete STA can at least converge to a \textit{local minimum}.\\
\indent \textit{Proof}. Let suppose the maximum number of iterations (denoted by \textit{M}) is big enough, then there must exist a number $N < M$, when $k > N$, no update of the better solution will happen. That is to say, $f(x_k) = f(x^{best}), \forall ~ k > N$, where $x^{best}$ is the solution in the \textit{N}th iteration. The $x^{best}$ is just the local minimum solution as $x^{*}$, and $|f(x_k) - f(x^{best})| = 0$.
\subsection{Global search ability of discrete STA}
It is not difficult to find that whether the discrete STA converges to a global minimum depends on the $x^{best}$ in the \textit{N}th iteration. If the $x^{best}$ is the global minimum, then according to the ``greedy criterion" we use to update the \textit{best} in pseudocode, when $k > N$, \textit{best} will always be $x^{best}$. In other words, the global convergence performance has much to do with the three operators we have designed.\\
\textbf{Theorem 2} the discrete STA can converge to a \textit{global minimum} in probability.\\
\indent \textit{Proof}. Let suppose $x^{*} = (a_1, \cdots, a_n)$ is the global minimum solution, and $x_k = (b_1, \cdots, b_n)$ is the \textit{k}th best solution. Then we discuss separately that how the state $x_k$ can transform to the best state $x^{*}$.\\
\indent If there exists the same sub-sequence in both $x_k$ and $x^{*}$, for instance, $x_k = (a_n, a_2, \cdots, a_{n-1},  a_1)$, then both swap transformation with small swap factor and appropriate shift transformation have the probability to swap or shift other positions in the sequence so that other positions in $x_k$ can transform exactly to the corresponding positions in the $x^{*}$.\\
\indent If there exists no same sub-sequence in $x_k$ and $x^{*}$, for instance, $x_k = (a_{n-1}, a_{n-2}, \cdots, a_2,  a_1)$, then swap transformation with big swap factor and proper symmetry transformation have the probability to swap or symmetrize all positions in the sequence so that all positions in $x_k$ can transform exactly to the corresponding positions in the $x^{*}$.
\subsection{Time complexity of discrete STA}
Because of the NP-hard property of TSP, it is impossible to solve the problem in polynomial time. As described in the section 1, to obtain a ``good solution" in a reasonable time is a optional choice, and this is the same case as the discrete STA, which aims at getting a satisfactory solution in as short a time as possible. In the pseudocode as described above, we can find that in the out loop, there are \textit{M} iterations, and in the inner loop, there exist three times of $SE$ transformations. It is not difficult to find that  the time complexity of the discrete STA is $O(M \cdot SE)$, that is to say, the discrete STA can achieve a global optimum in polynomial time in probability.
\section{Experimental results and discussion}
From the MATLAB scripts in section 3, we can find that as the $m_\alpha$ increases, the exchanging transformation matrix will become less efficient. The same phenomenon can be observed in other two transformations, and repetition of the transformation matrices appears frequently.\\
\indent For simplicity and efficiency, in current version of STA for TSP, parameters for three transformation operators are made as smaller as possible, that is to say, $m_a = 2, m_b = 1, m_c =0$ are specified in the paper, to avoid the trouble of parameters' adjusting.\\
\indent In order to test the performance of the proposed STA in discrete version, some benchmark symmetric traveling salesman problems from [3] are utilized for experiment, including \textit{ulysses16.tsp}, \textit{att48.tsp} and \textit{berlin52.tsp}.
In the same time, SA, ACO, which are recognized as distinguished algorithms for TSP are used for comparison with STA.
In SA, we set initial temperature at 5000, cooling rate at 0.97, and in ACO, $\alpha = 1, \beta =5, \rho = 0.9$, where,
$\alpha, \beta$ are used to control the relative weight of pheromone trail and heuristic value, and $\rho$ is the pheromone trail decay coefficient. In ACO and STA, the number of ants or the search enforcement is 20, and maximum iterations is fixed at 200. Considering that SA is not usually population based algorithm, the maximum iteration is extended especially for fairness. The threshold, or the total number of iterations in SA is set at 4000.\\
\indent Programs are run independently for 20 trails for each algorithms in MATLAB R2010b (version of 7.11.0.584) on Intel(R) Core(TM) i3-2310M CPU @2.10GHz under Window 7 environment, and comparison results for STA with SA and ACO are listed in Table 1. Some statistics and the run time are utilized to evaluate the performance of algorithms. The \textit{best} means the minimum of the results, the \textit{worst} indicates the maximum, and then it follows the \textit{mean}, \textit{st.dev.} (standard deviation). The run time is the average time used in 20 trails, which is measured in seconds.
\begin{table}[!htbp]
\caption{Results for benchmark test problems}
\footnotesize
\begin{tabular}{{lllll}}
\hhline
       Problems        & Performance &    SA      &     ACO     &   STA    \\ \hline
                       & best        &  73.9998   &   74.6287   & 73.9876  \\
                       & mean        &  74.4481   &   76.0864   & 74.0779  \\
 \textit{ulysses16.tsp}& worse       &  75.5391   &   78.7728   & 74.5939  \\
                       & st.dev.     &  0.4105    &   1.1062    & 0.1626   \\
                       & time(s)     &  2.9975    &   11.3038   & 1.2223   \\ \hline
                       & best        &  3.5266e4  &   3.7015e4  & 3.3724e4 \\
                       & mean        &  3.9667e4  &   3.8449e4  & 3.4872e4 \\
 \textit{att48.tsp}    & worse       &  4.5887e4  &   3.9801e4  & 3.6205e4 \\
                       & st.dev.     &  2.7453e3  &   862.4546  & 668.7553 \\
                       & time(s)     &  14.7605   &   102.4784  & 3.0462   \\ \hline
                       & best        &  8.1864e3  &   8.2404e3  & 7.5444e3 \\
                       & mean        &  8.9838e3  &   8.7776e3  & 8.2472e3 \\
 \textit{berlin52.tsp} & worse       &  9.5858e3  &   9.1513e3  & 8.6305e3 \\
                       & st.dev.     &  380.1004  &   267.1124  & 273.4509 \\
                       & time(s)     &  139.8399  &   118.0948  & 3.3438   \\
\hhline
\end{tabular}
\end{table}
\\
\indent As can be seen from the Table, STA outperforms SA and ACO in almost every performance index. To be more specific, the results of each benchmark problem are discussed separately in the following.\\
\indent \textit{ulysses16.tsp}: The \textit{best} of the results is obtained by STA, with the sequence of (7     6    14    13    12    16     1     3     2     4     8    15     5    11     9    10), which gets the best length of route at 73.9876, and the best route is plotted in Fig.1. To be more careful, we can find that the \textit{worst} solution gained by STA is even better than the \textit{best} of ACO, which indicates the strong search capability of STA. The \textit{st.dev.} of STA is almost approaching zero, and it shows that STA is also stable for this test problem. By the way, the STA consumes the half time of SA, and 1/9 time of ACO. The curves of the average fitness are illustrated following in Fig.2. The data created by SA are condensed to the same iterations(the same method applied to other two results obtained by SA). It is interesting to find that the ups and downs of the curve of SA in the early stage, because SA accepts a relatively worse solution by probability. Take the condensed data of SA into consideration, we can perceive that SA needs quite a long time for its steadily descending trend. However for ACO, after quickly to each a good fitness, it is trapped into stagnation point. But for STA, neither of the phenomena occurs, and it keeps decreasing before not a short time.
\begin{figure}[!htbp]
\includegraphics[width=8cm,height=6cm]{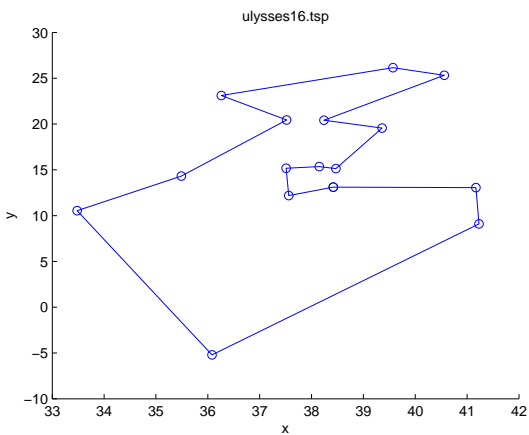}
\caption{the best route of ulysses16.tsp obtained by STA}
\end{figure}
\begin{figure}[!htbp]
\includegraphics[width=8cm,height=6cm]{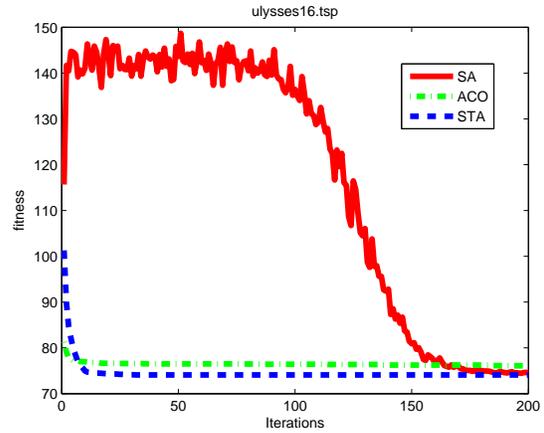}
\caption{curve of the average fitness for ulysses16.tsp}
\end{figure}
\\
\indent \textit{att48.tsp}: STA also acheives the \textit{best}, with the sequence of (9    40    15    12    11    23     3    22    16    41    34    48     5    29     2    42    26     4     35    45    10    24    32    39    25    14    13    21    47    20    33    46    36    30    43    17    27    19    37     6    28     7    18    44    31    38     8     1) and length of route at 3.3724e4, which is illustrated in Fig.3. The same situation is observed in the results between ACO and STA, that is, the \textit{worst} solution gained by STA is better than the \textit{best} of ACO. SA outperforms ACO in the \textit{best}, but \textit{mean} and \textit{worst} in ACO are better than that of SA. In this case, the \textit{st.dev.} of STA is not very satisfactory, although it is better than its counterparts. The time costed by STA is much shorter than other two algorithms, only around 1/4 and 1/30 time of that of SA and ACO consume, respectively. From Fig.4, we can find low degree of ups and downs in the curve of SA, and SA keeps good descending fitness in the first stage. ACO confronts the premature convergence again, few changes have happened in the late process. For STA, the fitness decrease sharply in the early stage, but can still keep decreasing in the later, which indicates the excellent performance of the designed operators in the discrete STA.
\begin{figure}[!htbp]
\centering
\includegraphics[width=8cm,height=6cm]{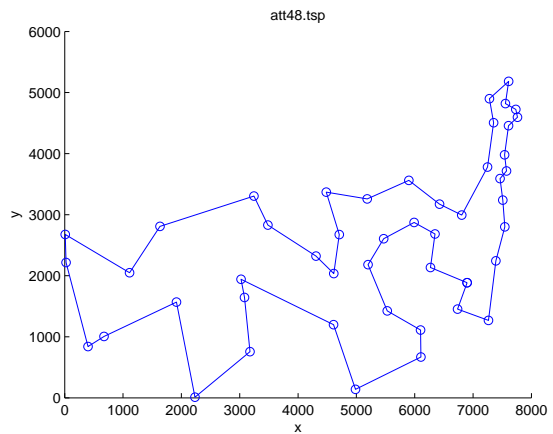}
\caption{the best route of att48.tsp obtained by STA}
\end{figure}
\begin{figure}[!htbp]
\centering
\includegraphics[width=8cm,height=6cm]{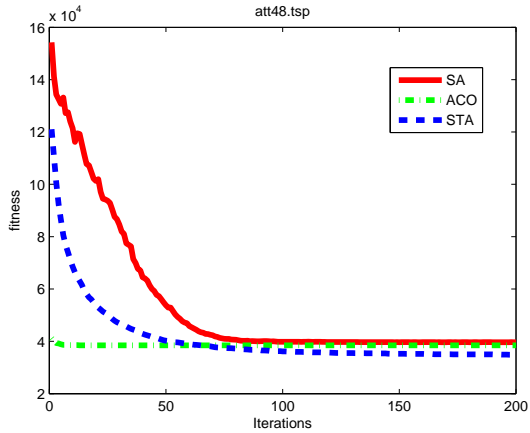}
\caption{curve of the average fitness for att48.tsp}
\end{figure}\\
\indent \textit{berlin52.tsp}: The STA wins the \textit{best} again, with the best sequence of (3    17    21    42     7     2    30    23    20    50    29    16    46    44    34    35    36    39    40    37    38    48    24     5    15     6     4    25    12    28    27    26    47    13    14    52    11    51    33    43    10     9     8    41    19    45    32    49     1    22    31    18) and length of route at 7.5444e3, which can be observed in Fig.5. At this time, ACO exhibits much better than SA in other performance except for the \textit{best}. But, STA achieves the best results on the whole, especially for the computational time. For the problem, STA consumes respectively 1/40 and 1/35 of the time costed by SA and ACO. As for TSP, the time complexity is really important, so the results gained by STA shows again that the discrete STA is really promising. In Fig.6, we can find that SA need quite a long time to reach a relatively good fitness, and then it becomes stagnated. On the contrary, the fitness of STA can decease quickly and keep descending till the end of the process.
\begin{figure}[!htbp]
\centering
\includegraphics[width=8cm,height=6cm]{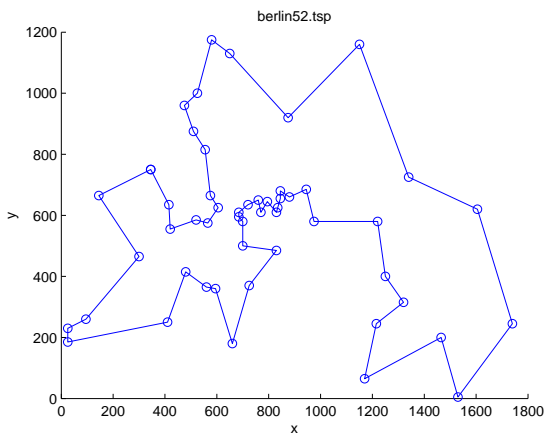}
\caption{the best route of berlin52.tsp obtained by STA}
\end{figure}
\begin{figure}[!htbp]
\centering
\includegraphics[width=8cm,height=6cm]{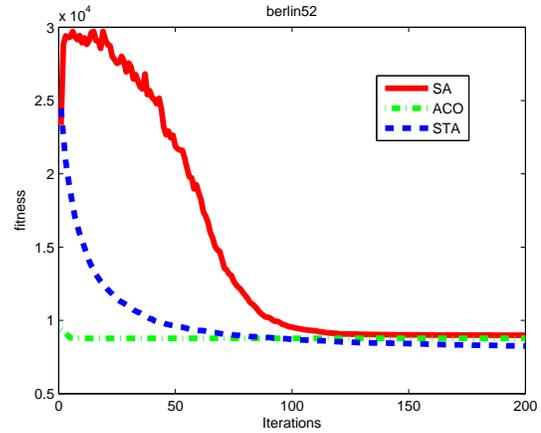}
\caption{curve of the average fitness for berlin52.tsp}
\end{figure}
\section{Conclusion}
Different from continuous search space, the space for traveling salesman problem is discrete, which corresponds to a
permutation of a sequence. In discrete version of state transition algorithm, three special state transition operators are
designed to manipulate the permutation. A simplest but much efficient parameters control strategy is proposed, in which, no parameters of the operators need adjusting. Some experiments are done to evaluate the proposed algorithm, and the results show that the discrete version of STA has much better performance not only in the search ability but also in the time consuming.\\
\indent Premature phenomenon is extensively existing in heuristic algorithms. To escape from a stagnation point, the paper focuses on the designing of operators. Accepting a relatively worse solution is a good idea, as seen in SA; however, it increases the
computing time and risks the non-convergence. It is really excited to see the fantastic performance of discrete STA only with ``greedy criterion", due to the excellent operators designed. On the other hand, difficulties are still to be confronted with large size problem, and the adjusting of the parameters need reconsidering. By the way, in current version of STA, three transformations are only for local permutation, and effective global permutation will be found in our future work as well as the equilibrium between them.

\end{document}